\documentclass[a4paper, 12pt]{article}
\usepackage{amsfonts}

 \usepackage{mathrsfs,amsfonts,amsmath,amssymb}
 \usepackage{color}

\makeatletter

\@addtoreset{equation}{section}
\makeatother
\newtheorem{Theorem}{Theorem}[section]
\newtheorem{Remark}[Theorem]{Remark}
\newtheorem{Lemma}[Theorem]{Lemma}

\newtheorem{Proposition}[Theorem]{Proposition}

\begin{document}

\title{\textbf{Stochastic continuity, irreducibility and non confluence
for SDEs with jumps}}
\author{Guangqiang Lan\footnote{Corresponding author. Supported by
China Scholarship Council, National Natural Science Foundation of
China (NSFC11026142) and Beijing Higher Education Young Elite Teacher Project (YETP0516).}
\\ \small School of Science, Beijing University of Chemical Technology, Beijing 100029, China
\\ \small Email: langq@mail.buct.edu.cn
\\  Jiang-Lun Wu
\\ \small Department of Mathematics, College of Science, Swansea University, Swansea SA2 8PP, UK
\\ \small Email: j.l.wu@swansea.ac.uk}
\date{}

\maketitle

\begin{abstract}
In this paper, we investigate stochastic continuity (with respect to the initial value),
irreducibility and non confluence property of the solutions of stochastic differential equations
with jumps. The conditions we posed are weaker than those relevant conditions existing in the
literature. We also provide an example to support our new conditions.
 \end{abstract}



\section{ Introduction and Main Results}

Given a probability space $(\Omega,\mathscr{F},P)$ endowed with a complete filtration $(\mathscr{F}_t)_{t\geq 0}$.
Let $d,m\in\mathbb{N}$ be arbitrarily fixed. We are concerned with the following stochastic differential equations (SDEs) with
jumps and with random coefficients
\begin{equation}\label{sdelan}\aligned X_t&=X_0+\int_0^t\sigma(s,\omega,X_s)dB_s+\int_0^tb(s,\omega,X_s)ds\\&
\quad+\int_0^{t+}\int_Uf_1(s,\omega,X_{s-},u)\tilde{N}_k(ds,du)+\int_0^{t+}\int_Uf_2(s,\omega,X_{s-},u)N_k(ds,du) \endaligned\end{equation}
where $B,\ N_k,\ \tilde{N}_k$ denote an $m$-dimensional $(\mathscr{F}_t)$-Brownian motion, a Poisson random measure and its compensated Poisson martingale measure, respectively, and $\mathbb{E}(N_k(ds,du))=ds\nu(du)$ with $\nu$ being a $\sigma$ finite measure on a given measurable space
$(U,\mathcal{B}(U))$, $\sigma:(t,\omega,x)\in[0,\infty)\times\Omega\times\mathbb{R}^d\mapsto\sigma(t,\omega,x)\in\mathbb{R}^d\otimes\mathbb{R}^m$ and $b:(t,\omega,x)\in[0,\infty)\times\Omega\times\mathbb{R}^d\mapsto b(t,\omega,x)
\in\mathbb{R}^d$ are progressively measurable functions, $f_i:(t,\omega,x,u)\in[0,\infty)\times\Omega\times\mathbb{R}^d
\times U\mapsto f_i(t,\omega,x,u)\in\mathbb{R}^d, i=1,2$ are $(\mathscr{F}_t)_{t\geq 0}$ predictable measurable functions
with
$$\mbox{supp} f_1(t,\omega,x,\cdot)\cap\mbox{supp} f_2(t,\omega,x,\cdot)
=\varnothing\, ,\quad \nu(\mbox{supp} f_2(t,\omega,x,\cdot))<\infty$$
for all $(t,\omega,x)\in[0,\infty)\times\Omega\times\mathbb{R}^d$,
and all the four functions are continuous with respect to the third variable $x$.

In order that the integrals in the definition of the solutions of the equation (\ref{sdelan}) are well-defined, we make
the following fundamental assumption which is enforced throughout the paper
\begin{equation}\label{t1}\mathbb{E}\int_0^T\sup_{|x|\le R}\Big[|b(s,\cdot,x)|+||\sigma(s,\cdot,x)||^2
+\int_U|f_i|^j(s,\cdot,x,u)\nu(du)\Big]ds<\infty
\end{equation}
for all $T,R>0, i,j=1,2$, where the norm $||\cdot||$ stands for the Hilbert-Schmidt norm $||\sigma||^2:=\sum\limits^d_{i=1}
\sum\limits^m_{j=1}\sigma^2_{ij}$ for any $d\times m$-matrix $\sigma=(\sigma_{ij})\in\mathbb{R}^d\otimes\mathbb{R}^m$
and $|\cdot|$ denotes the usual Euclidean norm on $\mathbb{R}^d$. As usual, we use $<\cdot,\cdot>$ to denote the
Euclidean inner product on $\mathbb{R}^d$.

Next, we fix $R>0$ arbitrarily. Let $\eta_R:\mathbb{R}_+\rightarrow\mathbb{R}_+$ be a differentiable function such that
$$\eta_R(0)=0,\ \eta'_R(x)\ge 0, \int_{0+}\frac{dx}{\eta_R(x)}=\infty.$$
We assume further that the coefficients of SDE (\ref{sdelan}) fulfill the following condition
\begin{equation}\label{t2}
\aligned&\quad||\sigma(t,\cdot,x)-\sigma(t,\cdot,y)||^2+2\langle x-y,b(t,\cdot,x)-b(t,\cdot,y)\rangle\\&
\quad+\sum_{i=1}^2\int_U|f_i(t,\cdot,x,u)-f_i(t,\cdot,y,u)|^2\nu(du)\\&
\quad+2\int_U\langle x-y,f_2(t,\cdot,x,u)-f_2(t,\cdot,y,u)\rangle\nu(du)
\\&\le g(t,\cdot)\eta_R(|x-y|^2),\ t\ge 0, |x|\vee|y|\le R, a.s.\endaligned
\end{equation}
for $g:[0,\infty)\times\Omega\to[0,\infty)$ to be a measurable function satisfying
$$\mathbb{E}\int_0^t g(s,\cdot)ds<\infty,\quad\forall t\ge0.$$

Under the above assumption (\ref{t2}), one can show (see e.g. \cite{Situ}) that there exists a unique pathwise solution
of SDE (\ref{sdelan}) which might blow up in finite time. In order to emphasize the solutions with different initial values,
we use the notation $X_t(x)$ for $t\ge0$ to denote the solution of SDE (\ref{sdelan}) starting from $X_0=x\in\mathbb{R}^d$.
Moreover, we denote the explosion time of the solution $X_t(x),t\ge0$, by
$$\zeta_x:=\inf\{t>0:|X_t(x)|=+\infty\}.$$

Our first main result concerns the uniformly stochastic continuity of the solution of SDE (\ref{sdelan}). We have the following

\begin{Theorem}\label{dingli1} Assume that the condition (\ref{t2}) holds. Let $X_t(x)$ and $X_t(y)$ be the solutions
of SDE (\ref{sdelan}) starting from $x,y\in\mathbb{R}^d$, respectively. Then for any $\varepsilon>0,$
$$\lim_{y\rightarrow x}P(\sup_{s< t\wedge\zeta_x\wedge\zeta_y}|X_s(x)-X_s(y)|\ge\varepsilon)=0,\quad \forall t>0.$$
If all solutions are non explosive, that is, $\zeta_x=\infty,\ \forall x\in \mathbb{R}^d$, then

$$\lim_{y\rightarrow x}P(\sup_{s\le t}|X_s(x)-X_s(y)|\ge\varepsilon)=0,\quad \forall t>0.$$

Namely, the solution of the SDE (\ref{sdelan}) is uniformly stochastic continuous with respect to the initial value
before any finite time $t$. 
\end{Theorem}

\begin{Remark}
Note that when the solutions are not global, the supremum must be taken with $s< t\wedge\zeta_x\wedge\zeta_y,$ otherwise it will 
be absurd since $|X_s(x)-X_s(y)|=\infty$ in this case.
\end{Remark}

Next, we consider the following SDE with deterministic coefficients (i.e., all coefficients are independent of $\omega\in\Omega$)
\begin{equation}\label{sdelan1}\aligned X_t&=x_0+\int_0^t\sigma(s,X_s)dB_s+\int_0^tb(s,X_s)ds\\&
\quad+\int_0^{t+}\int_Uf_1(s,X_{s-},u)\tilde{N}_k(ds,du)+\int_0^{t+}\int_Uf_2(s,X_{s-},u)N_k(ds,du). \endaligned
\end{equation}

Under the above condition (\ref{t2}), there exists a unique solution of equation (\ref{sdelan1}). Similar to Theorem 2.9.1 of \cite{Mao},
one can show that the solution is a Markov process. For $f\in C_b(\mathbb{R}^d)$, we define operator
$$P_{s,t}f(x):=\mathbb{E}^{s,x}f(X_t):=\int_{y\in\mathbb{R}^d} f(y)p_{s,t}(x,dy),\quad 0\le s\le t,\, x\in\mathbb{R}^d$$
where
$$p_{s,t}(x,A):=P(X_t\in A|X_s=x), \quad  0\le s\le t,\, x\in\mathbb{R}^d,\, A\in\mathcal{B}(\mathbb{R}^d)$$
is the transition probability measure of the Markov process. The operator family $\{P_{s,t}\}_{0\le s\le t}$ is the
Markov semigroup associated with the solution. Furthermore, our Theorem \ref{dingli1} ensures that the Markov
semigroup $\{P_{s,t}\}$ is Fellerian. That is, for any $f\in C_b(\mathbb{R}^d)$, we have $P_{s,t}f\in C_b(\mathbb{R}^d).$

We are going to study the irreducibility of the transition probability measure $p_{s,t}$ for any $0\le s\le t$. We say that the
family $\{p_{s,t}\}_{0\le s\le t}$ is {\it irreducible} if for any $0\le s\le t$ and $x\in \mathbb{R}^d$, $p_{s,t}(x,A)>0$
for any non empty open set $A\subset\mathbb{R}^d.$

Similar to \cite{zhang} and \cite{Qiao}, we introduce the following monotonicity and growth conditions. Suppose that
\begin{equation}\label{t9}
\aligned&\quad||\sigma(t,x)-\sigma(t,y)||^2+\langle x-y,b(t,x)-b(t,y)\rangle\\&
\quad+\sum_{i=1}^2\int_U|f_i(t,x,u)-f_i(t,y,u)|^2\nu(du)\\&
\quad+\int_U\langle x-y,f_2(t,x,u)-f_2(t,y,u)\rangle\nu(du)
\\&\le g(t)\eta(|x-y|^2),\ t\ge 0, a.s. \endaligned
\end{equation}
holds with $g\ge 0,\ \int_0^tg(s)ds<\infty,\forall t>0,$
$$\eta(x):=\left\{\begin{array}{ll} x\log\frac{1}{x},\qquad \qquad\qquad \qquad\quad\ \ x\le r<\frac{1}{e^2},\\
r\log\frac{1}{r}+(\log\frac{1}{r}-1)(x-r),\quad r<x,
\end{array}
\right.$$
\begin{equation}\label{t5}\aligned&||\sigma(t,x)||^2+2\langle x,b(t,x)\rangle
+\sum_{i=1}^2\int_U|f_i(t,x,u)|^2\nu(du)\\&
\quad+2\int_U\langle x,f_2(t,x,u)\rangle\nu(du) \leq
f(t)(|x|^2+1) \endaligned
\end{equation}
and
\begin{equation}\label{t6}|\sigma^T(t,x)x|^2+\sum_{i=1}^2\int_U(2\langle x,f_i(t,x,u)\rangle+|f_i(t,x,u)|^2)^2\nu(du)\leq
f(t)(|x|^2+1)^2
\end{equation}
holds for certain measurable function $f:[0,\infty)\to[0,\infty)$ with $\int_0^tf(s)ds<\infty,\forall t>0$.

It's obvious that there exists a unique non explosive solution of equation (\ref{sdelan1}) under conditions (\ref{t9}) and (\ref{t5}).
To investigate the irreducibility of $p_{s,t}$, we assume that
$m\ge d$, and  we need the following so called {\it strong
ellipticity condition} on the coefficient $\sigma$, that is, there exists $\lambda>0$ such that
\begin{equation}\label{t8}||\sigma^{-1}(t,x)||^2\le\lambda,\quad  t>0,x\in\mathbb{R}^d,
\end{equation}
where $\sigma^{-1}$ stands for the left inverse of matrix $\sigma.$

Our second main result is the following

\begin{Theorem}\label{dingli3}  Assume that the conditions (\ref{t9}), (\ref{t5}), (\ref{t6}) and (\ref{t8}) hold. If
there exists $2\le p<4$ such that
\begin{equation}\label{t7}||\sigma(t,x)||^2+\sum_{i=1}^2\int_U|f_i(t,x,u)|^2\nu(du)+(\int_U|f_2(t,x,u)|\nu(du))^2\leq
f(t)(|x|^p+1)
\end{equation}
hold with $f:[0,\infty)\to[0,\infty)$ being a measurable function satisfying $\int_0^t f^\frac{p}{2}(s)ds<\infty,\,\forall t\ge0$,
then $\{p_{s,t}\}$ is irreducible.
\end{Theorem}

\begin{Remark} \label{rm01}
It is worthwhile mentioning that here we do not need the assumption that $||\sigma(t,x)||$ is linear growth. The linear
growth condition on the coefficient $\sigma$ was required in both \cite{zhang,Qiao} while \cite{zhang} even only deals with
SDE without jumps. Our conditions (\ref{t5}), (\ref{t6}) and (\ref{t7}) are weaker than the linear
growth condition on the coefficient $\sigma$ (see Section 4) in \cite{zhang,Qiao}, even for relatively simpler SDEs without jumps in \cite{zhang}.
\end{Remark}

Our final task of the present paper concerns the non confluence property of the time-homogeneous SDE $(\ref{sdelan1})$
in which the coefficients are independent of $t$. We say that the solution $X_t$ of equation (\ref{sdelan1}) has {\it  non confluence property}, if for any initial values $x_0\neq y_0,$
$$P(X_t(x_0)\neq X_t(y_0),\ \forall t>0)=1.$$

In an early work \cite{Emery}, Emery studied such kind of non confluence property for general stochastic differential
equations without jumps under Lipschitzian coefficients. Yamada and Ogura considered in \cite{Yamada} for SDEs without
jumps with non-Lipschitz coefficients. We aim to give a new sufficient condition for the non confluence property of
the solution $X_t$ of the equation (\ref{sdelan1}).

Fix $R>0$ arbitrarily, let $\gamma_R:\mathbb{R}_+\rightarrow\mathbb{R}_+$ be a
differentiable function such that
$$\gamma_R(0)=0,\ \int_{0+}\frac{dx}{\gamma_R(x)}=\infty$$
and
$$\frac{x(\gamma'_R(x)+1)}{\gamma_R(x)}\le K,\ \forall x\in[0,\infty)$$
for some constant  $K> \frac{1}{2}$ which is independent of $x$ and $R$.

We have the following

\begin{Theorem}\label{dingli2} Assume that (\ref{t2}) holds with all the coefficients independent of $t$ and $\omega$. Let $K$ be given
as above. If for any $|x|\vee|y|\le R$
\begin{equation}\label{t3}
\aligned&||\sigma(x)-\sigma(y)||^2-\frac{2}{2K-1}\langle x-y,b(x)-b(y)\rangle
+\frac{1}{2K-1}\Big[\int_U|f_2(x,u)-f_2(y,u)|^2\\&-|f_1(x,u)-f_1(y,u)|^2+2\langle x-y,f_2(x,u)-f_2(y,u)\rangle\Big]\nu(du)
\leq \eta_R(|x-y|^2)\endaligned
\end{equation}
and
\begin{equation}\label{t4}|f_2(x,u)-f_2(y,u)|^2+2\langle x-y,f_2(x,u)-f_2(y,u)\rangle\ge0,
\end{equation}
then the unique solution of the time-homogeneous SDE $(\ref{sdelan1})$ has non confluence property.\end{Theorem}

\begin{Remark}
When there is no jumps, that is $f_i\equiv0, i=1,2$, in \cite{lan}, the authors showed that in one-dimensional case,
(\ref{t3}) also implies that the solution is stochastic monotonic. However, in the present case, there is no stochastic monotonicity
of the solution. Actually, we can conclude that if conditions (\ref{t2}) and (\ref{t3}) are satisfied, then
the process is stochastic monotonic between any two successive jumps.
\end{Remark}

The rest of the paper is organized as follows. In next section, Section 2, we show the uniformly
stochastic continuity with respect to initial value of the solution. Section 3 is devoted to
the proof of irreducibility of the transition probability $\{p_{s,t}\}_{0\le s\le t}$. In Section 4, we
present an example to illustrate that our conditions in Theorem \ref{dingli3} is indeed weaker than
those relevant known conditions in the literature. Finally in Section 5, we verify the non confluence property of solution
of the time-homogeneous equation (\ref{sdelan1}).

\section{Stochastic continuity with respect to initial value of the solution}

\textbf{Proof of Theorem \ref{dingli1}}

For any fixed $\varepsilon>0,$ let $x,y\in\mathbb{R}^d$ be such that $|y-x|<\varepsilon.$ Denote
$$\xi_t:=|\eta_t|^2:=|X_t(y)-X_t(x)|^2,\ \tau(x,y):=\inf\{t>0,\ \xi_t>\varepsilon^2\}.$$

Define the function
$\varphi_{\delta}:[0,\infty)\to[0,\infty)$ by
$$
\varphi_{\delta}(x):=\int_0^x\frac{ds}{\gamma_R(s)+\delta}.
$$

Then $\varphi{''}_{\delta}(x)\le 0, x>0.$ We can extend $\varphi_\delta$ to the real line (denoted by $\varphi_\delta$ again) such that
$\varphi{''}_{\delta}(x)\le 0, x\in\mathbb{R}.$
Denote
$$\tau_R(x,y):=\inf\{t,|X_t(x)|\vee|X_t(y)|>R\},$$
$$h_t:=b(t,\cdot,X_t(x))-b(t,\cdot,X_t(y)), e_t:=\sigma(t,\cdot,X_t(x))-\sigma(t,\cdot,X_t(y))$$
and
$$k_i(t-,u):=f_i(t,\cdot,X_{t-}(x),u)-f_i(t,\cdot,X_{t-}(y),u),i=1,2.$$

It's clear that $\tau_R(x,y)\rightarrow\zeta_x\wedge\zeta_y$ as $R\rightarrow\infty.$
By It\^o's formula, we have
$$\aligned\varphi_{\delta}(\xi_{t\wedge\tau(x,y)\wedge\tau_R})&=\varphi_\delta(|x-y|^2)+M_t+
\int_0^{t\wedge\tau(x,y)\wedge\tau_R(x,y)}\varphi'_{\delta}(\xi_s)[2\langle \eta_s,h_s\rangle+||e_s||^2\\&\quad+\int_U|k_1(s,u)|^2\nu(du)]ds
+2\int_0^{t\wedge\tau(x,y)\wedge\tau_R(x,y)}\varphi{''}_{\delta}(\xi_s)|e^T_s\eta_s|^2ds\\&
\quad+\int_0^{t\wedge\tau(x,y)\wedge\tau_R(x,y)}\int_U[\varphi_{\delta}(|\eta_s+k_2(s,u)|^2)-\varphi_{\delta}(\xi_s)]\nu(du)ds\\&
\quad+\int_0^{t\wedge\tau(x,y)\wedge\tau_R(x,y)}\int_U\Big[\varphi_{\delta}(|\eta_s+k_1(s,u)|^2)-\varphi_{\delta}(\xi_s)\\&\quad-
\varphi'_{\delta}(\xi_s)(|k_1(s,u)|^2+2\langle\eta_s,k_1(s,u)\rangle\Big]\nu(du)ds.
\endaligned$$
where $k_i(s,u)$ is defined similar to $k_i(s-,u)$ with $X_{s-}$ replaced by $X_s$.
Since $\varphi{''}_\delta\le0,$ by condition (\ref{t4}), we have
\begin{eqnarray*}&&\int_0^{t\wedge\tau(x,y)\wedge\tau_R(x,y)}\int_U[\varphi_{\delta}(|\eta_s+k_2(s,u)|^2)-\varphi_{\delta}(\xi_s)]\nu(du)ds\\
&\le&\int_0^{t\wedge\tau(x,y)\wedge\tau_R(x,y)}\int_U\varphi'_{\delta}(\xi_s)(|k_2(s,u)|^2+2\langle\eta_s,k_2(s,u)\rangle\nu(du)ds\end{eqnarray*}
and
$$\aligned&\int_0^{t\wedge\tau(x,y)\wedge\tau_R(x,y)}\int_U\Big[\varphi_{\delta}(|\eta_s+k_1(s,u)|^2)-\varphi_{\delta}(\xi_s)\\&\quad-
\varphi'_{\delta}(\xi_s)(|k_1(s,u)|^2+2\langle\eta_s,k_1(s,u)\rangle)\Big]\nu(du)ds\le0.\endaligned$$

Thus
$$\aligned
\mathbb{E}\varphi_{\delta}(\xi_{t\wedge\tau(x,y)\wedge\tau_R(x,y)})&\le\varphi_\delta(|x-y|^2)+\mathbb{E}\int_0^tg(s,\cdot)ds.
\endaligned$$

Taking
$\delta=|x-y|$ in the above inequality, we have
$$\mathbb{E}\varphi_{\delta}(\xi_{t\wedge\tau(x,y)\wedge\tau_R(x,y)})\leq \delta+\mathbb{E}\int_0^tg(s,\cdot)ds.$$

Hence
$$P(\tau(x,y)<t\wedge\tau_R(x,y))\varphi_{\delta}(\varepsilon^2)\leq\mathbb{E}\varphi_{\delta}(\xi_{t\wedge\tau(x,y)\wedge\tau_R(x,y)})\leq
\delta+C_t,$$
where $C_t=\mathbb{E}\int_0^tg(s,\cdot)ds.$ It follows that
$$P(\sup_{0\leq s<
t\wedge\tau_R(x,y)}|X_s(x)-X_s(y)|>\varepsilon)=P(\tau(x,y)<t\wedge\tau_R(x,y))\leq\varphi^{-1}_{\delta}(\varepsilon^2)(\delta+C_t).$$

Notice that the right hand side of the above inequality is independent of $R$. Let $R\rightarrow\infty$ and $\delta=|x-y|\rightarrow 0$ subsequently.
We then complete the proof. $\square$

\section{Irreducibility of $\{p_{s,t}\}$}\vskip0.2in

To investigate the irreducibility of $\{p_{s,t}\}$, we first introduce the following moment estimation of the maximal process.
\begin{Proposition}\label{yizhi} Assume (\ref{t5}) and (\ref{t6}) hold. Then for any $2\le p<4,$ the maximal process
$$Y_t:=\sup_{s\le t}|X_s|,\quad t\ge0$$
satisfies
$$\mathbb{E}({Y_t}^p)\le C_{t,p},\quad t\ge0.$$\end{Proposition}

To prove Proposition \ref{yizhi}, we need the following lemma.

\begin{Lemma}\label{martingale}
Suppose the assumptions of Proposition \ref{yizhi} hold. Let
$$M^c_t:=2\int_0^t\langle X_s,\sigma(s,X_s)dB_s\rangle,$$
$$M^d_t:=\sum_{i=1}^2\int_0^{t+}\int_U\Big(|f_i(s,X_{s-},u)|^2+2\langle X_{s-},f_i(s,X_{s-},u)\rangle\Big)\tilde{N}_k(ds,du)$$
and
$$M^{c*}_t:=\sup_{s\le t}|M^c_s|,\quad M^{d*}_t:=\sup_{s\le t}|M^d_s|.$$
Then for any $2\le p<4,$ there exist $K>0$ and $L>0$ such that
\begin{equation}\label{BDG}
\mathbb{E}((M^{c*}_t)^\frac{p}{2})\le C_p\Big(\frac{1}{2K}\mathbb{E}((Y^p_t+1))+\frac{K}{2}t^\frac{p-2}{2}
\int_0^tf^\frac{p}{2}(s)\mathbb{E}(Y^p_s+1)ds\Big)
\end{equation}
and
\begin{equation}\label{BDG1}
\mathbb{E}((M^{d*}_t)^\frac{p}{2})\le C{'}_p\Big(\frac{4-p}{2L^{4/(4-p)}}\mathbb{E}(Y^p_t+1)
+\frac{pL^{4/p}}{4}\int_0^tf(s)\mathbb{E}(Y^p_s+1)ds\Big).
\end{equation}

\end{Lemma}

\noindent \textbf{Proof} By Burkholder-Davis-Gundy inequality (for continuous martingales),
\[\aligned
\mathbb{E}((M^{c*}_t)^\frac{p}{2})&\le C_p\mathbb{E}[(\int_0^t|\sigma^T(s,X_s)X_s|^2ds)^\frac{p}{4}]\\&
\le C_p\mathbb{E}\Big((Y^2_t+1)^\frac{p}{4}\big(\int_0^t\frac{|\sigma^T(s,X_s)X_s|^2}{(|X_s|^2+1)}ds\big)^\frac{p}{4}\Big)\\&
\le C_p\Big(\frac{1}{2K}\mathbb{E}((Y^2_t+1)^\frac{p}{2})+\frac{K}{2}\mathbb{E}\big[\big(
\int_0^t\frac{|\sigma^T(s,X_s)X_s|^2}{(|X_s|^2+1)}ds\big)^\frac{p}{2}\big]\Big)\\&
\le C_p\Big(\frac{1}{2K}\mathbb{E}((Y^p_t+1))+\frac{K}{2}t^\frac{p-2}{2}
\int_0^tf^\frac{p}{2}(s)\mathbb{E}(Y^p_s+1)ds\Big).
\endaligned\]

We have used Young's inequality in the last second derivation and H\"older inequality in the last derivation.

On the other hand, by Burkholder-Davis-Gundy inequality for c\'adl\'ag martingales (see e.g., \cite{Qiao}), it follows that
$$\aligned\mathbb{E}((M^{d*}_t)^\frac{p}{2})&\le C{'}_p\sum_{i=1}^2\mathbb{E}\Big(\int_0^{t+}\int_UF_i(s-,u)^2N_k(ds,du)\Big)^\frac{p}{4}\\&
\le C{'}_p\sum_{i=1}^2\mathbb{E}\Big((Y^2_t+1)^{\frac{rp}{4}}(\int_0^{t+}\int_U\frac{F_i(s-,u)^2}{(|X_{s-}|^2+1)^r}N_k(ds,du))^\frac{p}{4}\Big)\endaligned$$
where
\begin{equation}\label{fi}F_i(s-,u):=|f_i(s,X_{s-},u)|^2+2\langle X_{s-},f_i(s,X_{s-},u)\rangle,\end{equation}
$r>0$ is a number to be determined later. By Young's inequality, we have
$$\aligned\mathbb{E}((M^{d*}_t)^\frac{p}{2})
\le C{'}_p\sum_{i=1}^2\mathbb{E}\Big(\frac{1}{aL^a}(Y^2_t+1)^{\frac{arp}{4}}+\frac{L^b}{b}
\Big(\int_0^{t+}\int_U\frac{F_i(s-,u)^2}{(|X_{s-}|^2+1)^r}N_k(ds,du)\Big)^\frac{bp}{4}\Big),\endaligned$$
where $a,b>0$, $\frac{1}{a}+\frac{1}{b}=1.$ Take $a,b,r$ such that
\begin{equation}\label{fangchengzu}\frac{1}{a}+\frac{1}{b}=1,\frac{bp}{4}=1\ \textrm{and}\ \frac{arp}{2}=p.\end{equation}

Then we have $b=\frac{4}{p}>1, a=\frac{4}{4-p}$ and $r=\frac{4-p}{2}.$ Thus, by condition (\ref{t6}),
it follows that

$$\aligned\mathbb{E}((M^{d*}_t)^\frac{q}{2})
&\le C'_p\sum_{i=1}^2\mathbb{E}\Big(\frac{4-p}{4L^{4/(4-p)}}(Y^p_t+1)+\frac{pL^{4/p}}{4}
\int_0^{t+}\int_U\frac{F_i(s-,u)^2}{(|X_{s-}|^2+1)^{(4-p)/2}}N_k(ds,du)\Big)\\&
\le C'_p\sum_{i=1}^2\mathbb{E}\Big(\frac{4-p}{4L^{4/(4-p)}}(Y^p_t+1)+\frac{pL^{4/p}}{4}
\int_0^{t}ds\int_U\frac{F_i(s,u)^2}{(|X_{s}|^2+1)^{(4-p)/2}}\nu(du)\Big)\\&
\le C'_p\Big(\frac{4-p}{2L^{4/(4-p)}}\mathbb{E}(Y^p_t+1)
+\frac{pL^{4/p}}{4}\int_0^tf(s)\mathbb{E}(Y^p_s+1)ds\Big).\endaligned$$

We complete the proof.  $\square$

\textbf{Proof of Proposition \ref{yizhi}}
By It\^o's formula, we have
\begin{equation}\label{ito} \aligned
|X_t|^2&=|x_0|^2+\int_0^t\Big(2\langle X_s,b(s,X_s)\rangle
+||\sigma(s,X_s)||^2\Big)ds\\&
\quad+\int_0^t\Big(\int_U(|f_1|^2(s,X_s,u)+F_2(s,u)\Big)\nu(du)ds\\&
\quad+2\int_0^t\langle X_s,\sigma(s,X_s)dB_s\rangle+\sum_{i=1}^2\int_0^{t+}\int_U F_i(s-,u)\tilde{N}_k(ds,du),
\endaligned\end{equation}
where $F_i(s-,u)$ is defined by (\ref{fi}), as in Lemma \ref{martingale}, and $F_i(s,u)$ is defined
in the same way with $X_{s-}$ replaced by $X_s$. Thus, by (\ref{t5}),
\begin{equation}\label{ito1} \aligned
Y_t^2\le|x_0|^2+\int_0^tf(s)(Y_s^2+1)ds+M^{c*}_t+M^{d*}_t.
\endaligned\end{equation}

Then we have
\begin{equation}\label{ito2} \aligned
\mathbb{E}(Y_t^p)&\le C_p\Big(|x_0|^p+(\int_0^tf(s)\mathbb{E}(Y_s^p+1)ds)
+\mathbb{E}((M^{c*}_t)^\frac{p}{2})+\mathbb{E}((M^{d*}_t)^\frac{p}{2})\Big)
\endaligned\end{equation}

By Lemma \ref{martingale} and (\ref{ito1}), we have
\begin{equation}\aligned
\mathbb{E}(Y_t^p)&\le C'_p\Big\{|x_0|^p+t^{\frac{p-2}{2}}(\int_0^tf^\frac{p}{2}(s)\mathbb{E}(Y_s^p+1)ds)\\&
\quad+\Big(\frac{1}{2K}\mathbb{E}(Y^p_t+1)+\frac{K}{2}t^\frac{p-2}{2}
\int_0^tf^\frac{p}{2}(s)\mathbb{E}(Y^p_s+1)ds\Big)\\&\quad
+\Big(\frac{4-p}{2L^{4/(4-p)}}\mathbb{E}(Y^p_t+1)
+\frac{pL^{4/p}}{4}\int_0^tf(s)\mathbb{E}(Y^p_s+1)ds\Big)\Big\}.
\endaligned\end{equation}

Set
$$\frac{C'_p}{2K}=\frac{(4-p)C'_p}{2L^{4/(4-p)}}=\frac{1}{4}.$$

We have $K=2C'_p, L=(C'_p(8-2p))^{(4-p)/4}.$ It follows that
\begin{equation}\aligned
\mathbb{E}(Y_t^p+1)&
\le A+B\int_0^t(f^\frac{p}{2}(s)+f(s))\mathbb{E}(Y^p_s+1)ds\Big\}
\endaligned\end{equation}
where $A=1+2C'_p|x_0|^p, B=C'_p((C'_p+1)t^{\frac{p-2}{2}}+\frac{p(C'_p(8-2p))^{(4-p)/p}}{4})$.
We then complete the proof by using Gronwall's lemma. $\square$

In what follows, we consider the irreducibility of $p_{s,t}$. For any $T>0,$ let us fix $t_1\in(0,T),$ whose value will be determined below. For any $\varepsilon>0$, define
$$X_{t_1}^\varepsilon:=X_{t_1}\cdot 1_{\{|X_{t_1}|\le\frac{1}{\varepsilon}\}}.$$

Then by Proposition \ref{yizhi}, for any $2\le p<4,$
$$\lim_{\varepsilon\downarrow0}\mathbb{E}|X_{t_1}^\varepsilon-X_{t_1}|^p=0.$$

For $t\in[t_1,T]$ and $y_0\in\mathbb{R}^d,$ define
$$Y_t^\varepsilon:=\frac{T-t}{T-t_1}X_{t_1}^\varepsilon+\frac{t-t_1}{T-t_1}y_0$$
and
$$h_t^\varepsilon:=\frac{y_0-X_{t_1}^\varepsilon}{T-t_1}-b(t,Y_t^\varepsilon).$$

Then
$$Y_{t_1}^\varepsilon=X_{t_1}^\varepsilon,\quad Y_{T}^\varepsilon=y_0$$
and
$$Y_{t}^\varepsilon=X_{t_1}^\varepsilon+\int_{t_1}^tb(s,Y_{s}^\varepsilon)ds+\int_{t_1}^th_{s}^\varepsilon ds,\quad t\in[t_1,T].$$

Consider the following SDE on $[t_1,T]$:
\begin{equation}\label{yt}\aligned Y_{t}&=X_{t_1}+\int_{t_1}^tb(s,Y_{s})ds+\int_{t_1}^th_{s}^\varepsilon ds
+\int_{t_1}^t\sigma(s,Y_{s})dB_s\\&\quad+\int_{t_1}^{t+}\int_Uf_1(s,Y_{s-},u)\tilde{N}_k(ds,du)
+\int_{t_1}^{t+}\int_Uf_2(s,Y_{s-},u)N_k(ds,du).\endaligned\end{equation}

We have the following

\begin{Proposition}\label{kongzhi} Suppose $b,\sigma$ and $f_i$ satisfy (\ref{t9}), (\ref{t5}), (\ref{t6}) and (\ref{t7}).
Then for any $T>0,$
\begin{equation}
\mathbb{E}|Y_T-y_0|^2\le C(t_1,T,p)^{e^{-2\int_{t_1}^T(g(s)+1)ds}}.
\end{equation}
\end{Proposition}

\textbf{Proof}\quad Set
$$Z_t^\varepsilon:=Y_{t}-Y_{t}^\varepsilon.$$

Since the coefficient $b$ is continuous with respect to $x$, and $h_t^\varepsilon$ is independent of $Y_t$ by definition, then conditions
(\ref{t9}) and (\ref{t5}) still hold when $b(t,x)$ is replaced by $b(t,x)+h_t^\varepsilon$. Thus, the SDE (\ref{yt}) has
a unique non explosive solution on $[t_1,T].$

By It\^o's formula and condition (\ref{t9}) we have
\begin{equation} \aligned\mathbb{E}|Z_t^\varepsilon|^2&=\mathbb{E}|X_{t_1}^\varepsilon-X_{t_1}|^2+
\mathbb{E}\Big\{\int_{t_1}^t\Big(2\langle Z_s^\varepsilon,b(s,Y_s)
-b(s,Y_s^\varepsilon)\rangle+||\sigma(s,Y_s)||^2\Big)ds\\&
\quad+\int_{t_1}^t\int_U(\sum_{i=1}^2\int_U|f_i(s,Y_s,u)|^2+2\langle Z_s^\varepsilon,f_2(s,Y_s,u)\rangle)\nu(du)ds\Big\}\\&
\le\mathbb{E}|X_{t_1}^\varepsilon-X_{t_1}|^2+\mathbb{E}\Big\{2\int_{t_1}^tg(s)\eta(|Z_s^\varepsilon|^2)ds
+2\int_{t_1}^t||\sigma(s,Y^\varepsilon_s)||^2ds\\&
\quad+2\int_{t_1}^t(\int_U(\sum_{i=1}^2\int_U|f_i(s,Y^\varepsilon_s,u)|^2
+|Z_s^\varepsilon||f_2(s,Y^\varepsilon_s,u)|\nu(du))ds\Big\}.
\endaligned\end{equation}

On the other hand, since
$$2\int_{t_1}^t\int_U|Z_s^\varepsilon||f_2(s,Y^\varepsilon_s,u)|\nu(du)ds\le\int_{t_1}^{t}|Z_s^\varepsilon|^2ds
+\int_{t_1}^t(\int_U|f_2(s,Y^\varepsilon_s,u)|\nu(du))^2ds,$$
and $\eta$ is concave and $\eta(x)\ge x$ for $r$ small enough by definition, by condition (\ref{t7}) and the fact that
$$|Y_s^\varepsilon|^p\le C_p(\sup_{s\le T}|X_s|^p+|y_0|^p),$$
we have
$$
\aligned\mathbb{E}|Z_t^\varepsilon|^2&\le\mathbb{E}|X_{t_1}^\varepsilon-X_{t_1}|^2+
\int_{t_1}^t[g(s)\eta(\mathbb{E}(|Z_s^\varepsilon|^2))+\mathbb{E}(|Z_s^\varepsilon|^2)]ds\\&
\quad+2(C_p(\mathbb{E}\sup_{s\le T}|X_s|^p+|y_0|^p)+1)\int_{t_1}^tf(s)ds\\&
\le C(t_1,T,p)+2\int_{t_1}^{t}(g(s)+1)\eta(\mathbb{E}(|Z_s^\varepsilon|^2))ds,\endaligned$$
where
\begin{equation}\label{changshu} C(t_1,T,p)=\mathbb{E}|X_{t_1}^\varepsilon-X_{t_1}|^2+2(C_p(\mathbb{E}
\sup_{s\le T}|X_s|^p+|y_0|^p)+1)\int_{t_1}^{T}f(s)ds.\end{equation}

Now, by utilising Bihari's inequality (see Lemma 2.1 of \cite{zhang}), we could ensure that
$$\mathbb{E}|Z_{t}^\varepsilon|^2\le C(t_1,T,p)^{e^{-2\int_{t_1}^T(g(s)+1)ds}}$$
holds for all $t\in[t_1,T].$ The proof is thus completed. \ $\square$

Now we are in a position to prove Theorem \ref{dingli3}.

\textbf{Proof of  Theorem \ref{dingli3}} Define
$$Y_t:=X_t,\quad t\in[0,t_1].$$

Then for any $t\in[0,T]$,
\begin{equation}\label{yt1}\aligned Y_{t}&=x_0+\int_{0}^tb(s,Y_{s})ds+\int_{0}^t1_{s>t_1}h_{s}^\varepsilon ds
+\int_{0}^t\sigma(s,Y_{s})dB_s\\&\quad+\int_{0}^{t+}\int_Uf_1(s,Y_{s-},u)\tilde{N}_k(ds,du)
+\int_{0}^{t+}\int_Uf_2(s,Y_{s-},u)N_k(ds,du).\endaligned\end{equation}

Define
$$\bar{B}_t:=B_t+\int_{0}^t1_{s>t_1}\sigma^{-1}(s,Y^\varepsilon_{s})h_{s}^\varepsilon ds$$
and
$$R_T^\varepsilon:=\exp\Big[-\int_{0}^T1_{s>t_1}\sigma^{-1}(s,Y^\varepsilon_{s})h_{s}^\varepsilon dB_s
-\frac{1}{2}\int_{0}^T1_{s>t_1}|\sigma^{-1}(s,Y^\varepsilon_{s})h_{s}^\varepsilon|^2 ds\Big].$$

Then by (\ref{t8}), the definition of $h_{s}^\varepsilon$ and the continuity of $b$ with respect to $x,$
we know that $R_T^\varepsilon\cdot P$ is a probability measure which is equivalent to $P$, and $\bar{B}_t$ is a $R_T^\varepsilon\cdot P$
Brownian motion. On the other hand, by \cite{Situ} Theorem 124 $\tilde{N}_k$ is still a Poisson martingale measure with the
same compensator $\nu(du)ds$ under the new probability measure $R_T^\varepsilon\cdot P$. By (\ref{yt1}), we have
\begin{equation}\label{yt2}\aligned Y_{t}&=x_0+\int_{0}^tb(s,Y_{s})ds
+\int_{0}^t\sigma(s,Y_{s})d\bar{B}_s\\&\quad+\int_{0}^{t+}\int_Uf_1(s,Y_{s-},u)\tilde{N}_k(ds,du)
+\int_{0}^{t+}\int_Uf_2(s,Y_{s-},u)N_k(ds,du).\endaligned\end{equation}

By the pathwise uniqueness of SDE (\ref{sdelan1}) (hence the uniqueness in law), $Y_\cdot(x_0)$ has the same
law as $X_\cdot(x_0)$ on $[0,T]$ for any $T>0$.
Thus we only need to prove that for each $0\le s< t, x\in \mathbb{R}^d$,
$$P(|Y_{t}(x_0)-Y_{s}(x_0)|>a)<1$$
for any $a>0$ since $R_T^\varepsilon\cdot P$ and $P$ are equivalent. Now
$$\aligned P(|Y_{t}(x_0)-Y_{s}(x_0)|>a)&\le P(|Y_{t}(x_0)-y|>\frac{a}{2})+P(|Y_{s}(x_0)-y|>\frac{a}{2})\\&
\le \frac{4}{a^2}(\mathbb{E}(|Y_{t}(x_0)-y|^2)+\mathbb{E}(|Y_{s}(x_0)-y|^2)).
\endaligned$$

According to Proposition \ref{kongzhi}, it follows that
$$\mathbb{E}(|Y_{t}(x_0)-y|^2)+\mathbb{E}(|Y_{s}(x_0)-y|^2)\le C(\tilde{t}_1,t,p)^{e^{-2\int_{\tilde{t}_1}^{t}(g(r)+1)dr}}
+C(\tilde{t}_2,s,p)^{e^{-2\int_{\tilde{t}_2}^{s}(g(r)+1)dr}}$$
where $\tilde{t}_1\le t,\tilde{t}_2\le s.$ Now let $\varepsilon$ to be sufficiently small,
$\tilde{t}_1$ close to $t$ and $\tilde{t}_2$ close to $s$. We have
$$P(|Y_{t}(x_0)-Y_{s}(x_0)|>a)<1.$$

This completes the proof.\ $\square$

\section{An example}

As pointed out in Remark \ref{rm01} in Section 1, our assumption on the coefficient $\sigma$ is weaker than those
relevant conditions carried out in \cite{zhang,Qiao}. Here let us give an example to support our conditions. We
create an example in the manner that it does satisfy our conditions (\ref{t5}), (\ref{t6}) and (\ref{t7}) but
it neither fulfill the condition ($H_2$) of Theorem 1.1 in \cite{zhang} nor the condition ($H_2$) of Theorem 1.3
in \cite{Qiao}. Thus our example indicates that our conditions are indeed weaker than those known conditions existing
in the literare.

\noindent\textbf{Example}\quad For simplicity, we only consider the time-homogeneous case with $f_1=f_2\equiv0$.
Suppose $d=m=2$. For any $2<p<4,$ define the $2\times2$-matrix coefficient $\sigma(x)$ and the drift vector coefficient
$b(x)$, respectively, by
\begin{equation}\sigma(x):=\left(\begin{array}{ccc}\frac{x_1}{1+|x|}&\frac{x_2}{1+|x|} \\ -(1+|x|^{\frac{p}{2}-1})x_2&(1+|x|^{\frac{p}{2}-1})x_1\end{array}\right)\end{equation}
and
\begin{equation}b(x):=-K_0(1+|x|^{p-2})x,\,\, \mbox{with constant}\,K_0\ge4.\end{equation}
Then, we have
\begin{eqnarray*}||\sigma(x)||^2+2\langle x,b(x)\rangle&=&\frac{|x|^2}{(1+|x|)^2}+(1+|x|^{\frac{p}{2}-1})^2|x|^2-2K_0(1+|x|^{p-2})|x|^2\\
&\le&1+2(1+|x|^{p-2})|x|^2-2K_0(1+|x|^{p-2})|x|^2\le1\\&\le& |x|^2+1\end{eqnarray*}
and
\begin{equation}||\sigma(x)||^2=\frac{|x|^2}{(1+|x|)^2}+(1+|x|^{\frac{p}{2}-1})^2|x|^2\le3(1+|x|^p).\end{equation}
On the other hand, it is clear that
$$\sigma(x)x=\left(\begin{array}{ccc}\frac{x_1}{1+|x|}&\frac{x_2}{1+|x|} \\ -(1+|x|^{\frac{p}{2}-1})x_2&(1+|x|^{\frac{p}{2}-1})x_1\end{array}\right)
\left(\begin{array}{ccc}x_1 \\ x_2\end{array}\right)=\left(\begin{array}{ccc}\frac{|x|^2}{1+|x|} \\ 0\end{array}\right).$$
So
\begin{equation}|\sigma(x)x|\le |x|\le 1+|x|^2.\end{equation}
Thus conditions (\ref{t5}), (\ref{t6}) and (\ref{t7}) hold in this case. We now show that condition (\ref{t9}) holds also.
We have indeed
$$\aligned||\sigma(x)-\sigma(y)||^2&=\sum_{i=1}^2(\frac{x_i}{1+|x|}-\frac{y_i}{1+|y|})^2
+\sum_{i=1}^2[(1+|x|^{\frac{p}{2}-1})x_i-(1+|y|^{\frac{p}{2}-1})y_i]^2\\&
\le|x-y|^2+2|x-y|^2+2(|x|^p+|y|^p-2(|x||y|)^{\frac{p}{2}-1}\langle x,y\rangle)\endaligned$$
and
$$\aligned\langle x-y,b(x)-b(y) \rangle&= K_0(2+|x|^{p-2}+|y|^{p-2})\langle x,y\rangle-K_0(|x|^2+|y|^2+|x|^{p}+|y|^{p})\\&
\le K_0(|x|^{p-2}+|y|^{p-2})|x||y|-K_0(|x|^{p}+|y|^{p}).\endaligned$$

Note that
$$(|x|^{p-2}+|y|^{p-2})|x||y|-(|x|^{p}+|y|^{p})=-(|x|-|y|)(|x|^{p-1}-|y|^{p-1})\le 0$$
and
$$\aligned &\quad(|x|^p+|y|^p-2(|x||y|)^{\frac{p}{2}-1}\langle x,y\rangle)-((|x|^{p}+|y|^{p})-(|x|^{p-2}+|y|^{p-2})\langle x,y\rangle)\\&
\le|x||y|(|x|^{\frac{p}{2}-1}-|y|^{\frac{p}{2}-1})^2.\endaligned$$

We have then
$$\aligned||\sigma(x)-\sigma(y)||^2+\langle x-y,b(x)-b(y)&\le3|x-y|^2-(K_0-2)(|x|-|y|)(|x|^{p-1}-|y|^{p-1})\\&
\quad+2|x||y|(|x|^{\frac{p}{2}-1}-|y|^{\frac{p}{2}-1})^2.\endaligned $$

On the other hand, if $|x|> |y|$, then
$$\aligned&\quad|x||y|(|x|^{\frac{p}{2}-1}-|y|^{\frac{p}{2}-1})^2-(|x|-|y|)(|x|^{p-1}-|y|^{p-1})\\&
\le|x||y|(|x|^{\frac{p}{2}-1}-|y|^{\frac{p}{2}-1})^2-|y|(|x|-|y|)(|x|^{p-2}-|y|^{p-2})\\&
=|y|(|x|^{\frac{p}{2}-1}-|y|^{\frac{p}{2}-1})[|x|(|x|^{\frac{p}{2}-1}-|y|^{\frac{p}{2}-1})-(|x|-|y|)(|x|^{\frac{p}{2}-1}+|y|^{\frac{p}{2}-1})]\\&
=|y|(|x|^{\frac{p}{2}-1}-|y|^{\frac{p}{2}-1})[|y|^{\frac{p}{2}-1}(|y|-|x|)+|x||y|(|x|^{\frac{p}{2}-2}-|y|^{\frac{p}{2}-2})]\le0\endaligned$$
where we have use the fact that $2<p<4$ in the last inequality. Next, by symmetry of $x$ and $y$, we know that
$$\aligned&\quad|x||y|(|x|^{\frac{p}{2}-1}-|y|^{\frac{p}{2}-1})^2-(|x|-|y|)(|x|^{p-1}-|y|^{p-1})\le0, \ \forall x,y\in\mathbb{R}^2.\endaligned$$
Hence
\begin{equation}\aligned||\sigma(x)-\sigma(y)||^2+\langle x-y,b(x)-b(y)&\le3|x-y|^2-(K_0-4)(|x|-|y|)(|x|^{p-1}-|y|^{p-1})\\&\le3|x-y|^2.\endaligned \end{equation}
We have thus verified the condition (\ref{t9}).

Since $\sigma\ge I,$ the condition (\ref{t8}) also holds.
By Theorem \ref{dingli3}, the transition probability $p_{s,t}$ is irreducible.
However, it is clear to see that there is no $K>0$ such that
$$||\sigma(x)||^2\le K(1+|x|^2)$$
and
$$|b(x)|\le K(|x|+1)$$
indicating that neither ($H_2$) of Theorem 1.1 in \cite{zhang} nor the condition ($H_2$) of Theorem 1.3
in \cite{Qiao} are fulfilled by our example.

\section{On the non confluence property of the solutions of the time-homogeneous SDEs}

Let us fix $x_0,y_0\in\mathbb{R}^d$ with $x_0\neq y_0$. For $0<\varepsilon<|x_0-y_0|$, we define
\begin{equation}\label{tingshi1}\hat{\tau}_\varepsilon:=\inf\{t>0,|X_t(x_0)-X_t(y_0)|\le\varepsilon\},\
\hat{\tau}:=\inf\{t>0,X_t(x_0)=X_t(y_0)\}. \end{equation}
It's obvious that
$\hat{\tau}_\varepsilon\rightarrow\hat{\tau}$, almost surely as $\varepsilon\rightarrow0.$

Next, we denote
\begin{equation}\label{tingshi2}\tau:=\inf\{t>0,|X_t(x_0)-X_t(y_0)|\ge 2|x_0-y_0|\}.
\end{equation}
Set the function
$$\varphi_{\delta}(x):=\exp\int_x^{c_0}\frac{ds}{\gamma_R(s)+\delta}.
$$
Then
$$\varphi'_\delta(x)=-\frac{\varphi_\delta(x)}{\gamma_R(x)+\delta}\le 0,\ \varphi{''}_\delta(x)
=\frac{\varphi_\delta(x)(1+\gamma'_R(x))}{(\gamma_R(x)+\delta)^2}.$$

By It\^o's formula, we have

$$\aligned\varphi_{\delta}(\xi_{t\wedge\tau\wedge\tau_R})&=\varphi_\delta(|x_0-y_0|^2)+M_t+
\int_0^{t\wedge\tau\wedge\tau_R}\varphi'_{\delta}(\xi_s)[2\langle \eta_s,h_s\rangle+||e_s||^2\\&\quad+\int_U|k_1(s,u)|^2\nu(du)]ds
+2\int_0^{t\wedge\tau\wedge\tau_R}\varphi{''}_{\delta}(\xi_s)|e^T_s\eta_s|^2ds\\&
\quad+\int_0^{t\wedge\tau\wedge\tau_R}\int_U[\varphi_{\delta}(|\eta_s+k_2(s-,u)|^2)-\varphi_{\delta}(\xi_s)]\nu(du)ds\\&
\quad+\int_0^{t\wedge\tau\wedge\tau_R}\int_U\Big[\varphi_{\delta}(|\eta_s+k_2(s,u)|^2)-\varphi_{\delta}(\xi_s)\\&\quad-
\varphi'_{\delta}(\xi_s)(|k_2(s,u)|^2+2\langle\eta_s,k_2(s,u))\Big]\nu(du)ds,
\endaligned$$
where $\tau_R=\inf\{t>0,|X_t(x_0)|\vee |X_t(y_0)|>R\}$.
By the definition of $\varphi_\delta$ and condition (\ref{t3}),
$$\aligned\varphi_{\delta}(\xi_{t\wedge\tau\wedge\tau_R})&
\le \varphi_\delta(|x_0-y_0|^2)+M_t+\int_0^{t\wedge\tau\wedge\tau_R}
\frac{\varphi_{\delta}(\xi_s)}{\gamma_R(\xi_s)+\delta}\Big[\frac{2(1+\gamma'_R(\xi_s))}{\gamma_R(\xi_s)+\delta}
|e^T_s\eta_s|^2\\&\quad-2\langle \eta_s,h_s\rangle-||e_s||^2+\int_U(|k_2(s,u)|^2-|k_1(s,u)|^2+2\langle\eta_s,k_2(s,u))\nu(du)\Big]ds\\&
\le \varphi_\delta(|x_0-y_0|^2)+(2K-1)\int_0^{t\wedge\tau\wedge\tau_R}
\frac{\varphi_{\delta}(\xi_s)}{\gamma_R(\xi_s)+\delta}\Big[||e_s||^2-\frac{2}{2K-1}
\langle \eta_s,h_s\rangle\\&\quad+\frac{1}{2K-1}\int_U(|k_2(s,u)|^2-|k_1(s,u)|^2+2\langle\eta_s,k_2(s,u))\nu(du)\Big]ds+M_t\\&
\le \varphi_\delta(|x_0-y_0|^2)+M_t+(2K-1)\int_0^t\varphi_{\delta}(\xi_s)ds,\endaligned$$
where
$$M_t=2\int_0^{t\wedge\tau\wedge\tau_R}\langle \eta_s,e_sdB_s\rangle+\sum_{i=1}^2\int_0^{{t\wedge\tau\wedge\tau_R}+}\int_U\Big(|k_i({s-},u)|^2+2\langle \eta_{s-},k_i({s-},u)\Big)\tilde{N}_k(ds,du)$$
is a real martingale. Take expectation on both sides. By Gronwall's lemma, we have
$$\mathbb{E}(\varphi_{\delta}(\xi_{t\wedge\tau\wedge\hat{\tau}_\varepsilon\wedge\tau_R}))\le\varphi_\delta(|x_0-y_0|^2)e^{(2K-1)t}.$$

On the other hand,
\begin{eqnarray*}
&& \mathbb{E}(\Phi_\delta(|X_{t\wedge\hat{\tau}_\varepsilon\wedge\tau\wedge\tau_R}(x_0)
-X_{t\wedge\hat{\tau}_\varepsilon\wedge\tau\wedge\tau_R}(y_0)|^2))\\
&\ge& \mathbb{E}(\Phi_\delta(|X_{t\wedge\hat{\tau}_\varepsilon\wedge\tau\wedge\tau_R}(x_0)
-X_{t\wedge\hat{\tau}_\varepsilon\wedge\tau\wedge\tau_R}(y_0)|^2)
1_{\hat{\tau}_\varepsilon\le
t\wedge\tau\wedge\tau_R})\\
&=&\Phi_\delta(\varepsilon^2)P(\hat{\tau}_\varepsilon\le
t\wedge\tau\wedge\tau_R).
\end{eqnarray*}

Thus,
\begin{equation}
P(\hat{\tau}_\varepsilon\le
t\wedge\tau\wedge\tau_R)\le C_t\exp{(-\int_{\varepsilon^2}^{\xi_0}\frac{ds}{\gamma(s)+\delta})},
\end{equation}
where the constant $C_t$ is independent of $R.$

Let $R\rightarrow\infty, \delta\rightarrow0, \varepsilon\rightarrow0 $ subsequently. We
have for any nonnegative $t,$
$$P(\hat{\tau}\le t\wedge\tau\wedge\zeta_{x_0}\wedge\zeta_{y_0})=0.$$

Let
$t\rightarrow\infty$, it follows that $P(\hat{\tau}\le \tau\wedge\zeta_{x_0}\wedge\zeta_{y_0})=0 $. Therefore, $\xi_\cdot$ is
positive almost surely on the interval $[0,\tau\wedge\zeta_{x_0}\wedge\zeta_{y_0}].$  Now we define
\begin{equation}\aligned &T_0:=0,\quad T_1:=\tau\wedge\zeta_{x_0}\wedge\zeta_{y_0},\\&
 T_2:=\inf\{t>T_1,|X_t(x_0)-X_t(y_0)|\le
|x_0-y_0|\}\wedge\zeta_{x_0}\wedge\zeta_{y_0},\endaligned\end{equation}
and generally
\begin{equation}\aligned &T_{2n}:=\inf\{t>T_{2n-1},|X_t(x_0)-X_t(y_0)|\le
|x_0-y_0|\}\wedge\zeta_{x_0}\wedge\zeta_{y_0},\\&
 T_{2n+1}:=\inf\{t>T_{2n},|X_t(x_0)-X_t(y_0)|\ge 2|x_0-y_0|\}\wedge\zeta_{x_0}\wedge\zeta_{y_0}.\endaligned\end{equation}

Due to Fang and Zhang \cite{Fang}, it is obvious that $T_n\rightarrow\zeta_{x_0}\wedge\zeta_{y_0}, a.s.$ as $n\rightarrow\infty$. Thus, $\xi_\cdot .0$ holds almost surely on 
$[T_{2n-1},T_{2n}] $ . By Theorem \ref{dingli1}, $X_t(x)$ is stochastic continuous with
 respect to the initial value $x$, thus the solution process $X_t(x)$ is a Feller process.
Further more, $\{X_t\}_{t\ge0}$ has the strong Markovian property since the process is right continuous
with left limit. Starting from $T_{2n} $ and applying the same arguments as in the first part of
the proof, $\xi_\cdot $ is also positive almost surely
on the interval $[T_{2n},T_{2n+1}]$. We complete the proof. $\square$

\textbf{Acknowledgement} The authors would like to thank Professor Feng-Yu Wang for useful discussions.

\end{document}